\documentclass[10pt]{amsart}

\usepackage{tikz}
\usepackage{amssymb,amsmath,latexsym,graphicx}
\usepackage{charter,eucal}
\usepackage{amssymb}
\usepackage[all,cmtip]{xy}
\usepackage{rotating,multicol}
\usetikzlibrary{decorations.markings}

\newtheorem{theorem}{Theorem }[section]

\newtheorem{lemma}[theorem]{Lemma}
\newtheorem{observation}[theorem]{Observation}

\newtheorem{remark}[theorem]{Remark}
\newtheorem{corollary}[theorem]{Corollary}
\newtheorem{proposition}[theorem]{Proposition}
\newtheorem{principle}[theorem]{\textsc{Principle}}

\newcommand{\bt}{\begin{theorem}}
\newcommand{\et}{\end{theorem}}
\newcommand{\bmt}{\begin{maintheorem}}
\newcommand{\emt}{\end{maintheorem}}
\newcommand{\bc}{\begin{corollary}}
\newcommand{\bl}{\begin{lemma}}
\newcommand{\ec}{\end{corollary}}
\newcommand{\el}{\end{lemma}}
\newcommand{\bo}{\begin{observation}}
\newcommand{\eo}{\end{observation}}
\newcommand{\bp}{\begin{proposition}}
\newcommand{\ep}{\end{proposition}}
\newcommand{\br}{\begin{remark}}
\newcommand{\er}{\end{remark}}
\newcommand{\bpr}{\begin{principle}}
\newcommand{\epr}{\end{principle}}

\def\PG{\mathbf{PG}}

\def\eop{\hspace*{\fill}$\blacksquare$}

\usepackage{amssymb} 
\usepackage{amsmath} 
\usepackage[utf8]{inputenc} 
\usepackage{latexsym}

\title {\bf On $k$-caps in $\PG(n, q)$, with $q$ even and $n \geq 3$}
\author {J. A. THAS \\ Ghent University}
\address{Ghent University, Department of Mathematics, Krijgslaan 281, S22, B-9000 Ghent, Belgium}
\email{jat@cage.ugent.be}
\begin {document}
\maketitle
\begin{abstract}
Let $m_2 (n, q) $ be the maximum size of $k$ for which there exists a $k$-cap in $\PG(n, q)$, and let $m^ \prime_2(n, q)$ be the second largest value of $k$ for which there exists a complete $k$-cap in $\PG(n, q)$. In this paper Chao's upper bound $q^2 - q + 5$ for $m^\prime _ 2(3, q) $, $q$ even and $q \ge 8$, will be improved. As a corollary new bounds for $m_2(n, q)$, $q$ even, $q \geq 8$ and $n \geq 4$, are obtained. Cao and Ou published a better bound but there seems to be a gap in their proof. \\

\end{abstract}

\section{Introduction}
A {\em k-arc } of $\PG(2, q)$ is a set of $k$ points, no three of which are collinear; a {\em k-cap} of $\PG(n, q)$, $n \ge 3$, is a set of $k$ points, no three of which are collinear. A $k$-arc or $k$-cap is {\em complete} if it is not contained in a ($k+1$)-arc or ($k+1$)-cap. The largest value of $k$ for which a $k$-arc of $\PG(2, q)$, or a $k$-cap of $\PG(n, q)$ with $n\ge 3$, exists is denoted by $m_2(n, q)$. The size of the second largest complete $k$-arc of $\PG(2, q)$ or $k$-cap of $\PG(n, q)$, $n\ge 3$, is denoted by $m^\prime_2(n, q)$. \\

\begin{theorem}
\label{thm1.1}
\begin {itemize}
\item [{\rm (i)}] $m_2(2, q) = q+2$, $q$ even  \cite {JWPH: 98};
\item [(ii)] $m_2(3, q) = q^2+1$, $q$ even, $q > 2$ \cite {JWPH: 85, RCB: 47, BQ: 52};
\item [(iii)] $m_2(n, 2) = 2^n$ \cite {RCB: 47};
\item [(iv)] $m_2(4, 4) = 41$ \cite {YE: 99};
\item [(v)] $m^\prime_2(n, 2) = 2^{n-1} + 2^{n-3}$ \cite {AAD: 90};
\item [(vi)] $m^\prime_2(3, 4) = 14$ \cite {JWPH: 87}.
\end {itemize}
\end{theorem}


\begin{theorem} [\cite {BS: 87, JAT: 87, JWPH: 98}]
\label{thm1.2}
Let $K$ be a $k$-arc of $\PG(2, q)$, $q$ even and $q > 2$, with $q - \sqrt{q} + 1 < k \le {q+1}$. Then $K$ can be uniquely extended to a ($q+2$)-arc of $\PG(2, q)$. 
\end{theorem}

For any   $k$-arc $K$ in $\PG(2, q)$ or $k$-cap $K$ in $\PG(n, q)$, $n\ge 3$, a {\em tangent} of $K$ is a line which has exactly one point in common with $K$. Let $t$ be the number of tangents of $K$ through a point $P$ of $K$ and let $\sigma_1(Q)$ be the number of tangents of $K$ through a point $Q \not\in K$. Then for a $k$-arc $K$ $t + k = q + 2$ and for a $k$-cap $K$ $t + k = q^{n-1} + q^{n-2} + \cdots + q + 2$. \\

\begin{theorem} [\cite {JWPH: 87}]
\label{thm1.3}
If $K$ is a complete $k$-arc in $\PG(2, q)$, $q$ even, or a complete $k$-cap in $\PG(n, q)$, $n\ge 3$ and $q$ even, then $\sigma_1(Q) \le t$ for each point $Q$ not on $K$.
\end{theorem}

\begin{theorem}[\cite {JMC: 99}]
\label{thm1.4}
\begin{equation}
\label{eq1}
m^\prime_2(3, q) \le q^2 - q + 5, q\ \text{even}, q \ge 8 .
\end{equation}
\end{theorem}

To prove Theorem \ref{thm1.4} J.-M. Chao relies on the following crucial lemma. \\

\begin{lemma} [\cite {JMC: 99}]
\label{lem1.5}
Let $K$ be a complete $k$-cap in $\PG(3, q)$ with $q$ even. If $\Pi$ is a plane such that $\vert \Pi \cap K\vert = x$, then
\begin{equation}
\label{eq2}
t(t - 1) \ge q(q + 2 - x)x.
\end{equation}
\end{lemma}

In the underlying paper the following improvement of Chao's result will be obtained.\\

\begin{theorem}[Main Theorem]
\label{thm1.6}
\begin{equation}
\label{eq3}
m^\prime_2(3, q) < q^2 - (\sqrt{5} - 1)q + 5, q\ \text{even}, q \ge 8.
\end{equation}
\end{theorem}

As a corollary new bounds for $m_2(n, q)$, $q$ even, $q \geq 8$ and $n \geq 4$, are obtained.

Combining the main theorem of \cite {SS: 93} with Theorem 1.6, there is an immediate improvement of the upper bound for $m^\prime_2(3, q)$, $q \geq 2048$. We thank T. Sz\H{o}nyi for bringing reference \cite {SS: 93} to our attention. 

\begin{theorem}
\label{thrm1.7}
\begin{equation}
\label{eq4}
m^\prime_2(3, q) < q^2 - 2q + 3\sqrt{q} + 2, q\ \text{even}, q \ge 2048.
\end{equation}
\end{theorem}

\section{A first improvement of Chao's bound}
\begin{theorem}
\label{thm2.1}
\begin{equation}
\label{eq5}
m^\prime_2(3, q) \le q^2 - q + 3, q\ \text{even}, q \ge 8.
\end{equation}
\end{theorem}
{\em Proof}. \quad
Let $K$ be a complete $k$-cap in $\PG(3, q)$, $q$ even, $q \ge 8$ and $k < q^2 + 1$.\\
Let $\Pi$ be a plane of $\PG(3, q)$ for which $4 \le \vert \Pi \cap K\vert \le q - 2 $.
Let $f(X) = q(q + 2 - X)X$. Then
\begin{equation}
\label{eq6}
t(t-1) \ge f(4) = f(q-2) = 4q(q - 2),
\end{equation}
by Lemma 1.5. So
\begin{equation}
\label{eq7}
t \ge \frac{1 + \sqrt{1+16q(q-2)}}{2} \ge 2q - \frac{7}{4} \ \text{for} \ q \ge 8.
\end{equation}
Hence $k \le q^2 + q + 2 - 2q + \frac{7}{4} = q^2 - q + \frac{15}{4}$, and consequently $k \le q^2 - q + 3$. \\

{\bf So we may assume that either $\vert \Pi \cap K \vert \le 3$, or $\vert \Pi \cap K \vert \ge q - 1$, for any plane $\Pi$ of $\PG(3, q)$. Let $l_1, l_2, ... , l_t$ be the $t$ tangents of $K$ through the point $P \in K$. We consider three cases depending on the number of planes containing $\l_i$ and intersecting $ K$ in at most 3 points}. \\

(\rm A) {\bf There exists exactly one plane $\Pi_{l_i}$ containing $l_i$ such that $\vert \Pi_{l_i }\cap K \vert \le  3, i = 1,2, ... , t$. We will show that in this case $k \leq q^2 - q + 3$.} \\

Assume there is exactly one plane $\Pi$ through $P$ with $\vert \Pi \cap K \vert \le 3$. Then for $i = 1, 2, ... , t$, $\Pi_{l_i} = \Pi$. Hence all tangents of $K$ containing $P$ are in $\Pi$. So $t \le q + 1$, a contradiction. Hence there are at least two planes $\Pi_1, \Pi_2$ through $P$ such that $\vert \Pi_i \cap K \vert \le 3, i= 1, 2$. Then $\vert \Pi_1 \cap \Pi_2 \cap K \vert = 2$. Consequently $t \ge 2(q - 1)$, and so $k \le q^2 + q + 2 - 2q + 2 = q^2 - q + 4$. 

Assume, by way of contradiction, that $k = q^2 - q + 4$. So $t = 2(q - 1)$. Then $\vert \Pi_1 \cap K \vert = \vert \Pi_2 \cap K \vert = 3$. All tangent lines at $P$ are contained in $\Pi_1$ and $\Pi_2$. Let $l$ be a tangent of $K$ at $P$ in $\Pi_1$, and consider the $q + 1$ planes containing $l$. The plane $\Pi_1$ is the only of these planes which intersects $K$ in 3 points, exactly $q - 1$ planes through $l$ contain 2 tangent lines at $P$ and so intersect $K$ in a $q$-arc and the remaining plane through $l$ contains exactly one tangent line at $P$ and so intersects $K$ in a $(q + 1)$-arc.

Let $\widetilde{\Pi}$ be the unique plane containing $l$ which intersects $K$ in a $(q + 1)$-arc, let $\widetilde{\Pi} \cap K = O$, and let $N$ be the kernel of $O$, that is, $N$ is the unique  point of $\widetilde{\Pi}$ which extends $O$ to a $(q + 2)$-arc of $\widetilde{\Pi}$. Clearly $N \in l$.

If $K^\prime$ is a $k^\prime$-arc of a plane $\PG(2, q)$ and $P^\prime \in \PG(2, q) \setminus K^\prime$, then the parity of the number of tangents of $K^\prime$ through $P^\prime$ is the parity of $k^\prime$, see Chapter 1 of \cite{JWPH: 98}. Hence, by considering $O$ and the $q - 1$ $q$-arcs whose planes contain $l$, we see that the number of tangents of $K$ through $N$ is at least $q + 1 + q - 1 = 2q$. As $K$ is complete we have $2q \le t$, so $k \le q^2 + q + 2 - 2q = q^2 - q + 2$, a contradiction.

Consequently $k \le q^2 - q + 3$.

(\rm B) {\bf Some  tangent $l_i, 1 \le i \le t$, is contained in at least two planes having at most three points in common with $K$}.\\

First we will prove that $k \leq q^2 - q + 5$. For $k = q^2 - q + 5$ and $k = q^2 -q + 4$ a contradiction will be obtained; the case $k = q^2 - q + 4$ will be subdivided in two cases. Hence it follows that also in Case (B) we have $k \leq q^2 - q + 3$. \\

Counting the points of $K$ on the $q+1$ planes containing $l_i$ gives
\begin{equation}
\label{eq8}
k - 1 \le 2.2 + (q - 1)q = q^2 - q + 4.
\end{equation}
So $k \le q^2 - q + 5$.

(\rm B.1) {\bf First, assume $k = q^2 - q + 5$}. Then two planes $\Pi_1, \Pi_2$ containing $l_i$ intersect $K$ in 3 points, while the remaining planes $\Pi_3, \Pi_4, ... , \Pi_{q + 1}$ containing $l_i$ intersect $K$ in $q + 1$ points. Let $l$ be a tangent of $K$ at $P$ in $\Pi_1$, distinct from $l_i$. Any plane $\zeta$ containing $l$, with $\zeta \not= \Pi_1$, intersects each $(q+1)$-arc $\Pi _i \cap K, i = 3, 4, ... , q + 1$, in exactly two points. Hence $\vert \zeta \cap K \vert \ge q$. Considering the lines $\zeta \cap \Pi_2$, we see that exactly two of the planes $\zeta$, say $\zeta _1$ and $\zeta_2$, intersect $K$ in $(q + 1)$-arcs $O_1$ and $O_2$, while the $q - 2$ other planes $\zeta$, say $\zeta_3, \zeta_4, ... , \zeta_q$, intersect $K$ in a $q$-arc.

Let $N_1$ be the kernel of $O_1$; then $N_1 \in l$. The number of tangents of $K$ containing $N_1$ is at least $q + 1 + q - 2 = 2q - 1$. As $K$ is complete we have $2q - 1 \le t$, so $k \le q^2 + q + 2 - 2q + 1 = q^2 - q + 3$, a contradiction.

(\rm B.2) {\bf Next, assume $k = q^2 - q + 4$}. Then, considering all planes containing $l_i$, there are two cases to consider.

(\rm B.2.1) {\bf Two planes $\Pi_1, \Pi_2$ containing $l_i$ intersect $K$ in three points, the plane $\Pi_3$ containing $l_i$ intersects K in $q$ points, and the remaining planes $\Pi_4, \Pi_5, ... , \Pi_{q + 1}$ containing $l_i$ intersect $K$ in $q + 1$ points.} 
Let $l$ be a tangent of $K$ at $P$ in $\Pi_1$, distinct from $l_i$. Any plane $\zeta$ containing $l$, distinct from $\Pi_1$, intersects each $(q + 1)$-arc $\Pi_i \cap K, i = 4, 5, ... , q + 1 $, in exactly two points; $q - 1$ of these planes $\zeta$ intersect $\Pi_3 \cap K$ in exactly two points. So for at least $q - 1$ of these planes $\zeta$ we have $\vert \zeta \cap K \vert \ge q$, and for all planes $\zeta$ we have $\vert \zeta \cap K \vert \ge {q - 1}$.

Assume that for all $q$ planes $\zeta$ we have $\vert \zeta \cap K \vert \ge q$. Let $s$ be the number of planes $\zeta$ for which $\vert \zeta \cap K \vert = q $ and let $u$ be the number of planes $\zeta$ for which $\vert \zeta \cap K \vert = q + 1$. Then
\begin{equation}
\label{eq9}
s(q - 1) + uq + 3 = q^2 - q + 4,\\
s + u = q.
\end{equation}
So $s(q - 1) + (q - s)q + 3 = q^2 - q + 4$, hence $s = q - 1$ and $u = 1$. Let $\zeta$ be the plane which intersects $K$ in a $(q + 1)$-arc $O$, and let $N \in l$ be the nucleus of $O$. The number of tangents of $K$ containing $N$ is at least $q + 1 + q - 1 =  2q$, so $k \le q^2 + q + 2 - 2q = q^2 - q + 2$, a contradiction.

So we may assume that for exactly $q - 1$ planes $\zeta$ we have $\vert \zeta \cap K \vert \ge q$ and that for exactly one plane $\zeta$ we have $\vert \zeta \cap K \vert = q - 1$. Assume that for $s$ planes $\zeta$ we have $\vert \zeta \cap K \vert = q$, and that for $u$ planes $\zeta$ we have $\vert \zeta \cap K \vert = q + 1$. Then
\begin{equation}
\label{eq10}
s(q - 1) + uq + q - 2 + 3 = q^2 - q + 4,
s + u = q -1.
\end{equation}
So $s(q - 1) + (q - 1 - s)q + q + 1 = q^2 - q + 4$, hence $s = q - 3$ and $u = 2$. Let $\zeta_1, \zeta_2$ be the planes containing $l$ which intersect $K$ in $(q + 1)$-arcs $O_1, O_2$, let $N_1, N_2$ be the nuclei of $O_1, O_2$, and let $\Pi_1 \cap K = \lbrace {P, P_1, P_2}\rbrace$. Assume first that $N_1 \not \in
P_1P_2$. Then the number of tangents of $K$ containing $N_1$ is at least $q + 1 + q - 3 + 2 = 2q$, so $k \le q^2 - q + 2$ a contradiction. Similarly if $N_2 \not \in P_1P_2$. Hence we may assume that $N_1 = N_2 = P_1P_2 \cap l$. Then the number of tangents of $K$ through $N_1$ is at least $q + 1 + q + q - 3 = 3q - 2$, so $k \le q^2 + q + 2 - 3q + 2 = q^2 - 2q + 4$, again a contradiction.

(\rm B.2.2) {\bf One plane $\Pi_1$ containing $l_i$ intersects $K$ in three points, and one plane $\Pi_2$ containing $l_i$ intersects $K$ in two points}. 
Consequently the other $q - 1$ planes $\Pi_3, \Pi_4, ... , \Pi_{q + 1}$ containing $l_i$ intersect $K$ in $q + 1$ points. Let $l$ be a tangent of $K$ at $P$ in $\Pi_1$, distinct from $l_i$. Any plane $\zeta$ containing $l$, distinct from $\Pi_1$, intersects each $(q + 1)$-arc $\Pi_i \cap K$, with $i = 3, 4, ... , q + 1$, in exactly two points. As $k = q^2 - q + 4$ it easily follows that for $q - 1$ of these planes $\zeta$ we have $ \vert \zeta \cap K \vert = q$, while for the remaining plane $\zeta$ we have $\vert \zeta \cap K \vert = q + 1$.

Let $\widetilde{\zeta}$ be the plane containing $l$ which intersects $K$ in a $(q + 1)$-arc $O$, and let $N$ be the nucleus of $O$. The number of tangents of $K$ containing $N$ is at least $q + 1 + q - 1 = 2q$, so $k \le q^2 - q + 2$, again a contradiction.

(\rm C) {\bf Some tangent $l_i$, with $1 \le i \le t$, is contained in no plane having at most three points in common with $K$}. \\

First we will prove that $k \leq q^2 - q + 5$. A contradiction will be obtained for $k \in \lbrace q^2 - q + 5, q^2 - q + 4 \rbrace$; for $k = q^2 - q + 4$ two cases have to be considered. Hence again $k \leq q^2 - q + 3$. \\

Then $\vert \Pi_j \cap K \vert \ge q - 1$ for each plane $\Pi_j$ containing $l_i$, with $j = 1, 2, ... , q + 1$. The arc $\Pi_j \cap K$ of $\Pi_j$ can be completed to a $(q + 2)$-arc of $\Pi_j$; see Theorem \ref{thm1.2}. This $(q + 2)$-arc meets $l_i$ in points $P, P_j$. As there are $q + 1$ points $P_j$ and $\vert l_i \setminus \lbrace P \rbrace \vert = q$, two of the points $P_j$ coincide, say $P_1 = P_2$. The number of tangents of $K$ containing $P_1$ is at least $2(q - 2) + 1 = 2q - 3$, so $k \le q^2 - q + 5$. \\

Now we make some observations on $(q - 1)$-arcs of $\PG(2, q)$, $q$ even. Let $\widetilde{K}$ be any $(q - 1)$-arc of $\PG(2, q)$, $q$ even, and let $\widetilde{l}$ be a tangent of $\widetilde{K}$ at $\widetilde{P} \in \widetilde{K}$. Let $\widetilde{C}$ be the unique $(q + 2)$-arc which contains $\widetilde{K}$; see Theorem \ref{thm1.2}. Put $\widetilde{C} \cap \widetilde{l} = \lbrace \widetilde{P}, \widetilde{N} \rbrace$. Then it is easy to see that exactly $q - 2$ points of $\widetilde{l} \setminus \lbrace \widetilde{P}, \widetilde{N} \rbrace$ are on exactly three tangents of $ \widetilde{K}$, and that exactly one point $\widetilde{R}$ of $\widetilde{l} \setminus \lbrace \widetilde{P}, \widetilde{N} \rbrace $ is on exactly one tangent of $\widetilde{K}$; also, $\widetilde{R} = \widetilde{l} \cap \widetilde{N}^\prime\widetilde{N}^{\prime\prime}$, with $\lbrace \widetilde{N}, \widetilde{N}^\prime, \widetilde{N}^{\prime\prime} \rbrace \cup \widetilde{K} = \widetilde{C}$.

(\rm C.1) {\bf First, assume $k = q^2 - q + 5$}. Then $\Pi_1 \cap K$ and $\Pi_2 \cap K$ are $(q - 1)$-arcs of $\Pi_1$ and $\Pi_2$. Let $r$ be the number of $(q - 1)$-arcs $\Pi_j \cap K$, let $s$ be the number of $q$-arcs $\Pi_j \cap K$ and let $u$ be the number of $(q + 1)$-arcs $\Pi_j \cap K$. Then
\begin{equation}
\label{eq11}
r(q - 2) + s(q - 1) + uq + 1 = q^2 - q + 5,
r + s + u = q + 1, \ \text{with}\ r \ge 2.
\end{equation}

So $r(q - 2) + s(q - 1) + (q + 1 - r - s)q + 1 = q^2 - q + 5$, hence $2r + s = 2q - 4$, with $r \ge 2$. If $s \ge 1$, then we have an extra tangent of $K$ containing $P_1$, so $k \le q^2 - q + 4$, a contradiction. Hence $s = 0, r = q - 2, u = 3$.


As the number of tangents of $K$ containing $P_1$ is exactly $2q - 3$, the nuclei of the three $(q + 1)$-arcs $\Pi_j \cap K$ are distinct from $P_1$. Let $N$ be one of these nuclei. Also, $P_1$ is on exactly one tangent of each of the $q - 4$ $(q - 1)$-arcs $\Pi_j \cap K$, distinct from the $(q - 1)$-arcs $\Pi_1 \cap K, \Pi_2 \cap K$. So $N$ is on at least three tangents of each of these $q - 4$ $(q - 1)$-arcs $\Pi_j \cap K$. Hence the number of tangents of $K$ containing $N$ is at least $2(q - 4) + q + 1 = 3q - 7 > 2q - 3$, a contradiction.

(\rm C.2) {\bf Finally, assume that $k = q^2 - q + 4$}. We have to consider two cases depending of the sizes of $\Pi_1 \cap K$ and $\Pi_2 \cap K$.

(\rm C.2.1) {\bf First, assume that $\Pi_1 \cap K$ and $\Pi_2 \cap K$ are $(q - 1)$-arcs}. The tangents of $K$ containing $P_1$ are the tangents of $\Pi_1 \cap K$ and $\Pi_2 \cap K$ containing $P_1$, and one extra tangent $l$. Assume that $l$ is a tangent of $\Pi_3 \cap K$. If $\Pi_3 \cap K$ is a $(q + 1)$-arc O, then $P_1$ is the nucleus of O, so there arise $q$ extra tangents, a contradiction; if $\Pi_3 \cap K$ is a $(q - 1)$-arc $K^\prime$, then $P_1$ is contained in at least three tangents of $K^\prime$, again a contradiction. Hence $\Pi_3 \cap K$ is a $q$-arc. Also, $\Pi_j \cap K$, with $j = 4, 5, ... , q + 1$, cannot be a $q$-arc. Let $r$ be the number of $(q - 1)$-arcs $\Pi_j \cap K$, and let $u$ be the number of $(q + 1)$-arcs $\Pi_j \cap K$. Then
\begin{equation}
\label{eq12}
r(q - 2) + uq + q - 1 + 1 = q^2 - q + 4,
r + u + 1 = q + 1.
\end{equation}
So $r(q - 2) + (q - r)q + q = q^2 - q + 4$, hence $r = q - 2$ and $u = 2$. Let $O_1, O_2$ be the $(q + 1)$-arcs $\Pi_j \cap K$, and let $N_1, N_2$ be the nuclei of $O_1, O_2$. Then $N_i \neq P_1, i = 1, 2$. Also $P_1$ is contained in exactly one tangent of each of the $q - 4$ $(q - 1)$-arcs $\Pi_j \cap K$, with $j \neq1, 2$. Hence the number of tangents of $K$ containing $N_1$ is at least $2(q -4) + q + 1 = 3q - 7 > 2q - 2$, clearly a contradiction.

(\rm C.2.2) {\bf Consequently, we may assume that $\Pi_1 \cap K$ is a $(q - 1)$-arc and that $\Pi_2 \cap K$ is a $q$-arc}
. Let $r$ be the number of $(q - 1)$-arcs $\Pi_j \cap K$, let $s$ be the number of $q$-arcs $\Pi_j \cap K$ and let $u$ be the number of $(q + 1)$-arcs $\Pi_j \cap K$. Then
\begin{equation}
\label{eq13}
r(q - 2) + s(q - 1) + uq + 1 = q^2 - q + 4,
r + s + u = q + 1, r \ge 1, s \ge 1.
\end{equation}
So $2r + s = 2q - 3, r \ge 1, s \ge 1$. Clearly, $s = 1$, as otherwise we have an extra tangent containing $P_1$, and then $k < q^2 - q + 4$. Hence $r = q - 2, s = 1, u = 2$. The nuclei of the two $(q + 1)$-arcs $\Pi_j \cap K$ are distinct from $P_1$. Let $N$ be one of these nuclei. Also, $P_1$ is on exactly one tangent of each of the $q - 3$ $(q - 1)$-arcs $\Pi_j \cap K$ distinct from $\Pi_1 \cap K$. So $N$ is on at least three tangents of each of these $q - 3$ $(q - 1)$-arcs $\Pi_j \cap K$. Consequently the number of tangents of $K$ containing $N$ is at least $2(q - 3) + q + 1 = 3q - 5 > 2q - 2$, a final contradiction. \eop \\

\section{Main Theorem}

\begin{theorem}
\label{thm 3.1}
\begin{equation}
\label{eq14}
m^\prime_2(3, q) < q^2 - (\sqrt{5} - 1)q + 5 , q\ \text{even}, q \ge 8.
\end{equation}
\begin{equation}
\label{eq15}
m^\prime_2(3, 4) = 14.
\end{equation}
\end{theorem}
{\em Proof} \quad
By [8] we have $m^\prime_2(3, 4) = 14$, and by Theorem 2.1 we have $m^\prime_2(3, 8) \leq 59$, which proves Theorem 3.1 for $q = 8$. So from now on we assume $q > 8$.

Let $K$ be a complete $k$-cap in $\PG(3, q)$, $q$ even, $q > 8$,  and $k < q^2 + 1$. Let $\Pi$ be a plane of $\PG(3, q)$ for which
\begin{equation}
\label{eq16}
5 \le \vert \Pi \cap K \vert \le q - 3.
\end{equation}
Let $f(X) = q(q + 2 - X)X$. Then by Lemma 1.5 of Chao
\begin{equation}
\label{eq17}
t(t - 1) \ge f(5) = f(q - 3) = 5q(q - 3).
\end{equation} 
So 
\begin{equation}
\label{eq18}
t \ge \frac{1 + \sqrt{1 + 20q(q -3)}}{2}.
\end{equation}
Put $\frac{1 + \sqrt{1 + 20q(q - 3)}}{2} \ge \sqrt{5}q - \alpha$, that is,
\begin{equation}
\label{eq19}
\sqrt{1 + 20q(q - 3)} \ge 2\sqrt{5}q - 2\alpha - 1.
\end{equation}
For $\alpha \le \sqrt{5}q - (1/2)$ this is equivalent to
\begin{equation}
\label{eq20}
1 + 20q(q - 3) \ge 20q^2 + 4\alpha^2 + 1 - 8\alpha \sqrt{5}q - 4\sqrt{5}q + 4\alpha,
\end{equation}
or
\begin{equation}
\label{eq21}
0 \ge 4\alpha^2 + \alpha(-8\sqrt{5}q + 4) + 60q - 4\sqrt{5}q,
\end{equation}
or
\begin{equation}
\label{eq22}
0 \ge \alpha^2 + \alpha(-2\sqrt{5}q + 1) + 15q - \sqrt{5}q.
\end{equation}
Put $\alpha = 3$. Then there arises $0 \ge 9 + 3(-2\sqrt{5}q + 1) + 15q - \sqrt{5}q$, that is, $0 \ge 12 + 15q - 7\sqrt{5}q$. This inequality is satisfied for $q > 16$.

Hence for $q > 16$ we have $t \ge \sqrt{5}q - 3$, and so,
\begin{equation}
\label{eq23}
k \le q^2 + q + 2 - \sqrt{5}q + 3,
\end{equation}
that is,
\begin{equation}
\label{eq24}
k \le q^2 + (1 - \sqrt{5})q + 5.
\end{equation}

For $q = 16$ it follows from (\ref{eq18}) that $t > 32$ and so $k \le 241$, which is equivalent to $k \le q^2 + (1 - \sqrt{5})q + 5$ with $q = 16$.

{\bf From now on suppose that either $\vert \Pi \cap K \vert \le 4$ or $\vert \Pi \cap K \vert \ge q - 2$ for any plane $\Pi$ of $\PG(3, q)$. Let $l_1, l_2, ... , l_t$ be the $t$ tangents of $K$ containing the point $P \in K$. Assume, by way of contradiction, that $k > q^2 + (1 - \sqrt{5})q + 5$. We consider three cases depending on the number of planes containing $l_i$ and intersecting $K$ in at most 4 points. In each case a contradiction will be obtained }.

(\rm A) {\bf Assume, by way of contradiction, that each $l_i$ is contained in exactly one plane $\Pi_{l_i}$ for which $\vert \Pi_{l_i} \cap K \vert \le 4$, with $i = 1, 2, ... , t$}. \\

(\rm A.1) {\bf Assume that there is exactly one plane $\Pi$ through $P$ with $\vert \Pi \cap K \vert \le 4$}. Then for $i = 1, 2, ... , t$ we have $\Pi_{l_i} = \Pi$. So $t \le q + 1$, hence $k \ge q^2 + 1$, a contradiction. \\

(\rm A.2) {\bf There are at least two planes $\Pi_1, \Pi_2$ through $P$ such that $\vert \Pi_i \cap K \vert \le 4, i = 1, 2$}. Then $\vert \Pi_1 \cap \Pi_2 \cap K \vert = 2$. Consequently $t \ge 2(q - 2)$, and so $k \le q^2 + q + 2 - 2q + 4 = q^2 - q + 6$.

The plane $\Pi_1$ intersects $K$ in a $m$-arc, $m \le 4$, and contains at least $q - 2$ tangents of $K$ at $P$. Let $P_1 \in (K \cap \Pi_1) \setminus P$ and assume that $PP_1$ is contained in $\alpha$ planes $\Pi$ with $\vert \Pi \cap K \vert \le 4$. Then $t \ge \alpha(q - 2)$, so $k \le q^2 + (1 - \alpha)q + 2 + 2\alpha$. Consequently
\begin{equation}
\label{eq25}
q^2 + (1 - \alpha)q + 2 + 2\alpha > q^2 + (1 - \sqrt{5})q + 5,
\end{equation}
or
\begin{equation}
\label{eq26}
(\sqrt{5} - \alpha)q + 2\alpha - 3 > 0
\end{equation}
This gives a contradiction for $\alpha > 2$ with $q > 8$. So $PP_1$ is contained in at most two planes intersecting $K$ in at most four points.

Assume, by way of contradiction, that for some plane $\Pi$ of $\PG(3, q)$ we have $\Pi \cap K = \lbrace P \rbrace$. As there are at least two planes $\Pi, \Pi^\prime$ through $P$ intersecting $K$ in at most four points, we have $\vert \Pi \cap \Pi^\prime \cap K \vert = 2$ and so $\vert \Pi \cap K \vert \ge 2$, a contradiction.

Let $\PG(2, q)$ be a plane of $\PG(3, q)$ not containing $P$ and let $\sigma$ be the projection of $\PG(3, q) \setminus \lbrace P \rbrace$ from $P$ onto $\PG(2, q)$. Further, let $\mathcal{P}$ be the set of all images under $\sigma$ of all points of $K \setminus \lbrace P \rbrace$ contained in planes $\Pi$, with $P \in \Pi$, for which $\vert \Pi \cap K \vert \le 4$, and let $\mathcal{B}$ be the set of all images under $\sigma$ of the sets $\Pi \setminus  \lbrace P \rbrace$. Then there arises an incidence structure $(\mathcal{P}, \mathcal{B})$ of points and lines for which
\begin{itemize}
 \item [(1)] $\vert \mathcal{B} \vert \ge 2$,
 \item [(2)] any two distinct lines in $\mathcal{B}$ have exactly one point in common,
 \item [(3)] each point is contained in at most two lines,
 \item [(4)] each line contains at most three points and at least one point.
 \end{itemize}
 It follows easily that $2 \le \vert \mathcal{B} \vert \le 4$. For each value of $\vert \beta \vert$ we will find a contradiction.

($\alpha$)  \underline {$\vert \mathcal{B} \vert = 4$} \\
Then $t = 4(q - 2)$, so $k = q^2 + q + 2 - 4q + 8 = q^2 - 3q + 10$. Hence $q^2 - 3q + 10 > q^2 + (1 - \sqrt{5})q + 5$, or $5 > (4 - \sqrt{5})q$, a contradiction as $q > 8$.

($\beta$) \underline {$\vert \mathcal{B} \vert = 3$} \\
If $\vert \mathcal{P} \vert = 3$, then $t = 3(q - 1)$, so $k = q^2 - 2q + 5$. Hence $q^2 - 2q + 5 > q^2 + (1 - \sqrt{5})q + 5$, or $(3 - \sqrt{5})q < 0$, a contradiction.

If $\vert \mathcal{P} \vert = 4$, then $t = 2(q - 1) + q - 2$, so $k = q^2 - 2q + 6$. Hence $q^2 + (1 - \sqrt{5})q + 5 < q^2 - 2q + 6$, or $(3 - \sqrt{5})q - 1 < 0$, a contradiction.

If $\vert \mathcal{P} \vert = 5$, then $t = q - 1 +2(q - 2)$, so $k = q^2 - 2q + 7$. Hence $q^2 + (1 - \sqrt{5})q + 5 < q^2 - 2q + 7$, or $(3 - \sqrt{5})q < 2$, a contradiction.

If $\vert \mathcal{P} \vert = 6$, then $t = 3(q - 2)$, so $k = q^2 - 2q + 8$. Hence $q^2 + (1 - \sqrt{5})q + 5 < q^2 - 2q + 8$, or $(3 - \sqrt{5})q < 3$, a contradiction.

($\gamma$) \underline {$\vert \mathcal{B} \vert = 2$} \\
By Theorem 2.1 we may assume that $k \le q^2 - q + 3$.

If $\vert \mathcal{P} \vert = 1$, then $t = 2q$, so $k = q^2 - q + 2$.

If $\vert \mathcal{P} \vert = 2$, then $t = 2q - 1$, so $k = q^2 - q + 3$.

If $\vert \mathcal{P} \vert = 3$, then $t = 2q - 2$, so $k = q^2 - q + 4$, a contradiction.

If $\vert \mathcal{P} \vert = 4$, then $t = 2q - 3$, so $k = q^2 - q + 5$, a contradiction.

If $\vert \mathcal{P} \vert = 5$, then $t = 2q - 4$, so $k = q^2 - q + 6$, a contradiction.

Hence the cases $k = q^2 - q + 2$ and $k = q^2 - q + 3$ have still to be considered.

($\gamma .1$) \underline {$k = q^2 - q + 2$} \\
On $K$ there are two points $P, P_1$ such that $PP_1$ is contained in two planes $\Pi_1, \Pi_2$ intersecting $K$ in just $\lbrace P, P_1 \rbrace$, and in $q - 1$ planes $\Pi_3, \Pi_4, ... , \Pi_{q + 1}$ intersecting $K$ in a $(q + 2)$-arc.

Let $P^\prime \in (\Pi_3 \cap K) \setminus \lbrace P, P_1 \rbrace$ and let $l$ be a tangent of $K$ at $P^\prime$. Assume, by way of contradiction, that each plane containing $l$ intersects $K$ in a $m$-arc with $m > 4$, so $m \ge q - 2$. These $m$-arcs $K^\prime_i$, with $i = 1, 2, ... , q + 1$, are extendable to $(q + 2)$-arcs $C_i$. Let $ C_i \cap l = \lbrace N_i, P^\prime \rbrace, i = 1, 2, ... , q + 1$. At least two of the points $N_1, N_2, ... , N_{q + 1}$ coincide, say $N_1 = N_2$. A plane $\Pi^\prime$ containing $l$, but not containing $P$ nor $P_1$, intersects each of the $(q + 2)$-arcs $\Pi_i \cap K$, with $i = 3, 4, ... , q+ 1$, in either 0 or 2 points; so $\vert \Pi^\prime \cap K \vert$ is even. A plane $\Pi^\prime$ containing $l$ and either $P$ or $P_1$ intersects $K$ in $q$ points. Hence each plane containing $l$ intersects $K$ in a $m$-arc, with $m$ even. Counting tangents of $K$ containing $N_1$, we obtain at least $2(q - 3) + 1 + q - 1 = 3q - 6$ tangents. So $k \le q^2 + q + 2 - 3q + 6 = q^2 - 2q + 8$, a contradiction for $q > 8$. We conclude that there is a plane $\Pi^\prime$ containing $l$ with $\vert \Pi^\prime \cap K \vert \le 4$.

Assume, by way of contradiction, that $l$ is contained in at least two planes $\Pi^\prime, \Pi^{\prime \prime}$ with $\vert \Pi^\prime \cap K \vert \le 4, \vert \Pi^{\prime \prime} \cap K \vert \le 4$. Then, by a previous argument, these intersections have an even number of points and so $\vert \Pi^\prime \cap K \vert \in \lbrace 2, 4 \rbrace$ and $\vert \Pi^{\prime \prime} \cap K \vert \in \lbrace 2, 4 \rbrace$. Now we count the points of $K$ in planes containing $l$, and obtain $k \le (q - 1)(q - 1) + 7 = q^2 - 2q + 8$, a contradiction for $q > 8$.

Hence $l$ is contained in exactly one plane $\Pi^\prime$ for which $\vert \Pi^\prime \cap K \vert \le 4$. It follows that the roles of $P$ and $P^\prime$ may be interchanged.

Let $l^\prime$ be a second tangent of $K$ containing $P^\prime$, with $l^\prime \not \subset \Pi^\prime$. Let $\tilde K = K \cap \Pi_3, \Pi^\prime \cap \tilde K = \lbrace P^\prime, P^\prime_1 \rbrace$. If $P^\prime_1 \in \lbrace P, P_1 \rbrace$, then $\vert \Pi^\prime \cap K \vert = q$, a contradiction. Hence $P^\prime_1 \not \in \lbrace P, P_1 \rbrace$. With $P^\prime$ there corresponds an incidence structure $(\mathcal{P}^\prime, \mathcal{B}^\prime)$ of points and lines. As $k = q^2 - q + 2$, we necessarily have $\vert \mathcal{P}^\prime \vert = 1$ and $\vert \mathcal{B}^\prime \vert = 2$. Hence $\Pi^\prime \cap K = \lbrace P^\prime, P^\prime_1 \rbrace$. If $ \tilde \Pi^\prime$ is the unique plane containing $l^\prime$ and intersecting $K$ in at most 4 points, then $\tilde \Pi^\prime \cap K = \lbrace P^\prime, P^\prime_1 \rbrace$. Also, the roles of $P$ and $P^\prime_1$, $P^\prime$ and $P^\prime_1$, $P$ and $P_1$ can be interchanged.

Interchanging $\Pi_3$ and $\Pi_i$, $i \in \lbrace 3, 4, .... , q + 1 \rbrace$, and interchanging $P^\prime$ with any point of $(\Pi_i \cap K) \setminus \lbrace P, P_1 \rbrace$, we see that $K$ is partitioned into $(q^2 - q + 2)/2$ pairs, where each pair is contained in two planes intersecting $K$ in that pair and in $q - 1$ planes intersecting $K$ in a $(q + 2)$-arc. Any other plane contains either 0 or $q$ points of $K$. Each point $Q$ of $K$ is contained in $2q$ tangents; the two planes on $Q$ intersecting $K$ in two points each contain $q$ of these tangents.

Now we count the planes intersecting $K$ in a $(q + 2)$-arc, and obtain
\begin{equation}
\label{eq27}
\frac{q^2 - q + 2}{2} .(q - 1) / \frac{q + 2}{2}.
\end{equation}
Hence $q + 2  \vert  (q^2 - q + 2)(q - 1)$, so $q + 2  \vert  24$, that is $q \in \lbrace 2, 4 \rbrace$, a contradiction.

($\gamma.2$) \underline {$k = q^2 - q + 3$} \\
Then on $K$ there are points $P, P_1$ such that $PP_1$ is contained in two planes $\Pi_1, \Pi_2$ with $\Pi_1 \cap K = \lbrace P, P_1 \rbrace, \Pi_2 \cap K = \lbrace P, P_1, P_2 \rbrace$, and in $q - 1$ planes $\Pi_3, \Pi_4, ... , \Pi_{q + 1}$ intersecting $K$ in a $(q + 2)$-arc.

Let $P^\prime \in ( \Pi_3 \cap K) \setminus \lbrace P, P_1 \rbrace$ and let $l$ be a tangent of $K$ at $P^\prime$. Assume, by way of contradiction, that each plane containing $l$ intersects $K$ in a $m$-arc with $m > 4$, so $m \ge q - 2$. These $m$-arcs $K^\prime_i$, with $i = 1, 2, ... , q + 1$, are extendable to $(q + 2)$-arcs $C_i$. Let $C_i \cap l = \lbrace N_i, P^\prime \rbrace, i = 1, 2, ... , q + 1$. At least two of the points $N_1, N_2, ... , N_{q + 1}$ coincide, say $N_1 = N_2$. A plane $\Pi^\prime$ containing $l$, but not containing $P$ nor $P_1$, intersects each of the $(q + 2)$-arcs $\Pi_i \cap K$, with $i = 3, 4, ... , q + 1$, in either 0 or 2 points. So if $P_2 \not \in \Pi^\prime$, then $\vert \Pi^\prime \cap K \vert$ is even. A plane $\Pi^\prime$ containing $l$ and either $P$ or $P_1$, but not $P_2$, intersects $K$ in $q$ points. Hence $q$ planes containing $l$ intersect $K$ in a $m$-arc, with $m$ even. Counting tangents of $K$ containing $N_1$, we obtain at least $2(q - 3) + 1 + q - 2 = 3q - 7$ tangents. So $k \le q^2 + q + 2 - 3q + 7 = q^2 - 2q + 9$, a contradiction for $q > 8$. We conclude that there is a plane $\Pi^\prime$ containing $l$ with $\vert \Pi^\prime \cap K \vert \le 4$.

Assume, by way of contradiction, that $l$ is contained in at least two planes $\Pi^\prime, \Pi^{\prime \prime}$ with $\vert \Pi^\prime \cap K \vert \le 4, \vert \Pi^{\prime \prime} \cap K \vert \le 4$. Now we count the points of $K$ in planes containing $l$, and obtain $k \le q^2 -2q + 9$, a contradiction for $q > 8$.

Hence $l$ is contained in exactly one plane $\Pi^\prime$ for which $\vert \Pi^\prime \cap K \vert \le 4 $. As all tangents of $K$ at $P_1$ are contained in $\Pi_1 \cup \Pi_2$, it follows that each tangent of $K$ at $P_1$ is contained in exactly one plane intersecting $K$  in at most 4 points. Hence all points of $K \setminus \lbrace P_2 \rbrace$ play the same role.

Let $l^\prime$ be a second tangent of $K$ containing $P^\prime$, with $l^\prime \not \subset \Pi^\prime$. Let $K \cap \Pi_3 = \tilde K, \Pi^\prime \cap \tilde K = \lbrace P^\prime, P^\prime_1 \rbrace$. If $P^\prime_1 \in \lbrace P, P_1 \rbrace$, then $\vert \Pi^\prime \cap K \vert \ge q$, a contradiction. Hence $P^\prime_1 \not \in \lbrace P, P_1 \rbrace$. With $P^\prime$ there corresponds an incidence structure $(\mathcal{P}^\prime, \mathcal{B}^\prime)$ of points and lines (see first part of (A)).

As $k = q^2 - q + 3$, we necessarily have $\vert \mathcal{P}^\prime \vert = 2$ and $\vert \mathcal{B}^\prime \vert = 2$. Hence $\vert \Pi^\prime \cap K \vert \in \lbrace 2, 3 \rbrace$ and $\Pi^\prime \cap K \supset \lbrace P^\prime, P^\prime_1 \rbrace$. Let $\widetilde \Pi^\prime$ be the unique plane containing $l^\prime$ and intersecting $K$ in at most 4 points, and let $\widetilde \Pi^\prime \cap \widetilde K = \lbrace P^\prime, \widetilde P^\prime_1 \rbrace$. If $P^\prime_1 \neq \widetilde P^\prime_1$, then by the structure of $(\mathcal{P}^\prime, \mathcal{B}^\prime)$ we have  $\lbrace P^\prime_1, \widetilde P^\prime_1 \rbrace \subset \Pi^\prime$, clearly a contradiction. Hence $P^\prime_1 = \widetilde P^\prime_1$, and so $\lbrace P^\prime, P^\prime_1 \rbrace \subset \widetilde \Pi^\prime \cap K$.

Without loss of generality we may assume that $\Pi^\prime \cap K = \lbrace P^\prime, P^\prime_1, P^\prime_2 \rbrace$ and $\widetilde \Pi^\prime \cap K = \lbrace P^\prime, P^\prime_1 \rbrace$. As $\vert \Pi^\prime \cap K \vert $ is odd, the set $\Pi^\prime \cap K$ has to contain the point $P_2$. Consequently $P_2 = P^\prime_2$.

Interchanging $\Pi_3$ and $\Pi_i$, $i \in \lbrace 3, 4, \cdots, q + 1 \rbrace$, and interchanging $P^\prime$ with any point of $(\Pi_i \cap K) \setminus \lbrace P, P_1 \rbrace$, we see that $K \setminus \lbrace P_2 \rbrace$ is partitioned into $(q^2 - q + 2)/2$ pairs, where each pair is contained in one plane intersecting $K$ in that pair, in one plane intersecting $K$ in that pair together with $P_2$, and in $q - 1$ planes intersecting $K$ in a $(q + 2)$-arc. Any other plane contains 0, 1, $q$ or $q + 1$ points of $K$.

Now we count the planes intersecting $K$ in a $(q + 2)$-arc and obtain
\begin{equation}
\label{eq28}
\frac{q^2 - q + 2}{2}.(q - 1)/\frac{q + 2}{2}.
\end{equation}
Hence $q + 2 \vert (q^2 - q + 2)(q - 1)$, so $q + 2 \vert 24$, that is $q \in \lbrace 2, 4 \rbrace$, a final contradiction.

We  conclude that there is some tangent $l_i$ containing $P$, with $i \in \lbrace 1, 2, ... , t \rbrace$, which is contained in exactly $\theta > 1$ planes having at most 4 points in common with $K$.

(\rm B) {\bf Assume, by way of contradiction, that some tangent $l$ of $K$ is contained in no plane intersecting $K$ in at most 4 points}.\\

Hence each plane $\Pi_i$ containing $l$ satisfies $\vert \Pi_i \cap K \vert \geq q - 2$, with $i = 1, 2, ... , q + 1$. By Theorem 1.2 the arc $\Pi_i \cap K$ can be extended to a $(q + 2)$-arc $C_i$; let $C_i \cap l = \lbrace N_i, P \rbrace$ with $l \cap K = \lbrace P \rbrace$. For at least two planes $\Pi_i$, say $\Pi_1$ and $\Pi_2$, we have $N_1 = N_2$.

(\rm B.1) {\bf First we prove that $N_1$ is on a tangent of $K$ not in $\Pi_1 \cup \Pi_2$; clearly $N_1$ is on at least $2q - 5$ tangents of $K$ contained in $\Pi_1 \cup \Pi_2$}. Assume the contrary. Then for any plane $\Pi_i \not \in \lbrace \Pi_1, \Pi_2 \rbrace$, the arc $\Pi_i \cap K$ must have an odd number of points. So $\Pi_i \cap K$ either is a $(q - 1)$-arc or a $(q + 1)$-arc, $i \in \lbrace 3, 4, ... , q + 1 \rbrace$. Also, $N_i \neq N_1$ for $i = 3, 4, ... , q - 1$. If $\Pi_i \cap K$ is a $(q - 1)$-arc and $C_i \setminus (\Pi_i \cap K) = \lbrace N_i, N^\prime_i, N^{\prime \prime}_i \rbrace, i \in \lbrace 3, 4, ... , q + 1 \rbrace$, then $N_1 \in N^\prime_i N^{\prime \prime}_i$, as otherwise $N_1N^\prime_i$ and $N_1N^{\prime \prime}_i$ are tangents of $\Pi_i \cap K$.

Let $r$ be the number of planes $\Pi_i$, with $i \neq 1, 2$, for which $\Pi_i  \cap K$ is a $(q - 1)$-arc, and let $s$ be the number of planes $\Pi_i$, with $i \neq 1, 2$, for which $\Pi_i \cap K$ is a $(q + 1)$-arc. The number of points of $K$ is at least
\begin{equation}
\label{eq29}
r(q - 2) + sq + 2(q - 3) + 1, {\rm with} \  r + s = q - 1.
\end{equation}
As $K$ is complete, by Theorem 2.1
\begin{equation}
\label{eq30}
r(q -2) + (q - 1 - r)q + 2(q - 3) + 1 \leq q^2 - q + 3,
\end{equation}
so
\begin{equation}
\label{eq31}
r \geq q - 4.
\end{equation}
We may assume that $\Pi_3 \cap K$ is a $(q - 1)$-arc. The number of tangents of $K$ containing $N_3$ is at least
\begin{equation}
\label{eq32}
q - 1 + 2(r - 1) \geq q - 1 + 2q - 10 = 3q - 11.
\end{equation}
Hence
\begin{equation}
\label{eq33}
k \leq q^2 + q + 2 - 3q + 11 = q^2 - 2q + 13.
\end{equation}
So
\begin{equation}
\label{eq34}
q^2 - 2q + 13 > q^2 + (1 - \sqrt{5})q + 5,
\end{equation}
a contradiction for $q > 8$.

Consequently $N_1$ is on a tangent $l^\prime$ of $K$ not in $\Pi_1 \cup \Pi_2$.

(\rm B.2) {\bf Now we consider all planes $\Pi^\prime_i$ containing the tangent $l^\prime$, with $i = 1, 2, ... , q + 1$. We will show that:\\
(a) For each plane $\Pi^\prime_i$ such that $\vert \Pi^\prime _i \cap K \vert \geq q - 2$ the point $N_1$ does not extend the arc $\Pi^\prime_i \cap K $. \\
(b) For each $i$ we have $\vert \Pi^\prime_i \cap K \vert \geq q - 2$}.\\

(a) Let $\vert \Pi_1 \cap K \vert = \alpha, q - 2 \leq \alpha \leq q + 1, \vert \Pi_2 \cap K \vert = \beta, q - 2 \leq \beta \leq q + 1$. Then $N_1$ is contained in at least $\alpha + \beta$ tangents of $K$. Now we consider all planes $\Pi^\prime_i$ containing the tangent $l^\prime$, with $i = 1, 2, ... , q + 1$. Assume, by way of contradiction, that $m = \vert \Pi^\prime_i \cap K \vert \geq q - 2$ and that the $(q + 2)$-arc $C^\prime_i$ extending $\Pi^\prime_i \cap K$ intersects $l^\prime$ in $\lbrace N_1, P^\prime \rbrace$, with $l^\prime \cap K = \lbrace P^\prime \rbrace, i \in \lbrace 1, 2, ... , q + 1 \rbrace$. Then the number of tangents of $K$ containing $N_1$ is at least
\begin{equation}
\label{eq35}
\alpha + \beta + m - 3 \geq 2q - 4 + m - 3 \geq 3q - 9.
\end{equation}
Hence
\begin{equation}
\label{eq36}
k \leq q^2 + q + 2 - 3q + 9 = q^2  - 2q + 11.
\end{equation}
So $q^2 - 2q + 11 > q^2 + ( 1 - \sqrt{5})q + 5$, a contradiction. Consequently for $\vert \Pi^\prime_i \cap K \vert \geq q - 2$ we have $N_1 \not \in C^\prime_i, i \in \lbrace 1, 2, ... , q + 1 \rbrace$. 

(b) Next, assume by way of contradiction that for at least one plane $\Pi^\prime_i$ containing $l^\prime$, say $\Pi^\prime_1$, we have $\vert \Pi^\prime_1 \cap K \vert \leq 4$. Let $\Pi^\prime_2$ be the plane $ll^\prime$. Now we count the points of $K$ in the planes $\Pi^\prime_i$, with $i = 1, 2, ... , q + 1$. Let \\
$\theta_1$ be the number of planes $\Pi^\prime_i, i \in \lbrace 3, 4, ... , q + 1 \rbrace$, containing a tangent of $\Pi_1 \cap K$ through $N_1$ and a tangent of $\Pi_2 \cap K$ through $N_1$,\\
$\theta_2$ be the number of planes $\Pi^\prime_i, i \in \lbrace 3, 4, ... , q + 1 \rbrace$, containing a tangent of $\Pi_1 \cap K$ through $N_1$, but no tangent of $\Pi_2 \cap K$ through $N_1$,\\
$\theta_3$ be the number of planes $\Pi^\prime_i, i \in \lbrace 3, 4, ... , q + 1 \rbrace$, containing a tangent of $\Pi_2 \cap K$ through $N_1$, but no tangent of $\Pi_1 \cap K$ through $N_1$,\\
$\theta_4$ be the number of planes $\Pi^\prime_i, i \in \lbrace 3, 4, ... , q + 1 \rbrace $, containing no one of the tangents of $\Pi_1 \cap K$ or $\Pi_2 \cap K$ through $N_1$.\\
Then, as $N_1 \notin C^\prime_i$ for $\vert \Pi^\prime_i \cap K \vert \geq q - 2$, we have
\begin{equation}
\label{37}
k \leq 4 + q - 1 + \theta_1(q - 2) + \theta_2(q - 1) + \theta_3(q - 1) + \theta_4q, 
{\rm with} \ 2 + \theta_1 + \theta_2 + \theta_3 + \theta_4 = q + 1.
\end{equation} 

Hence
\begin{equation}
\label{eq38}
k \leq q(\theta_1 + \theta_2 + \theta_3 + \theta_4) - (2\theta_1 + \theta_2 + \theta_3) + q + 3,
\end{equation}
so
\begin{equation}
\label{eq39}
k \leq q(q - 1) - (2\theta_1 + \theta_2 + \theta_3) + q + 3.
\end{equation}
Now we have \\
$\theta_1 + \theta_2 \geq \vert \Pi_1 \cap K \vert - 2 \geq q - 4,$ \\
$\theta_1 + \theta_3 \geq \vert \Pi_2 \cap K \vert - 2 \geq q - 4$.\\
Hence
\begin{equation}
\label{eq40}
k \leq q(q -1) - 2q + 8 + q + 3 = q^2 - 2q + 11.
\end{equation}
So $q^2 - 2q + 11 > q^2 + (1 - \sqrt5)q + 5$, a contradiction.

Hence no plane $\Pi^\prime_i$ containing $l^\prime$ intersects $K$ in a $m$-arc, with $m \leq 4, 1 \leq i \leq q + 1$. Consequently, for each plane $\Pi^\prime_i$ containing $l^\prime$ we have $\vert \Pi^\prime_i \cap K \vert \geq q - 2$. Also, we know that the $(q + 2)$-arc $C^\prime_i$ extending $\Pi^\prime_i \cap K$ does not contain $N_1$, with $i = 1, 2, ... , q + 1$.

(\rm B.3) {\bf A final contradiction will be obtained by considering the possible intersections $ \Pi^\prime_i  \cap K, i = 1, 2, ... , q + 1$}. It is easy to see that at least $q - 6$ planes $\Pi^\prime_i$ containing $l^\prime$ intersect $K$ in a $m$-arc having at least 3 tangents containing $N_1$; these planes are the planes containing $l^\prime$ passing through distinct tangents of $\Pi_1 \cap K$ and $\Pi_2 \cap K$ containing $N_1$. For any such plane $\Pi^\prime_i$ the arc $\Pi^\prime_i \cap K$ is either a $(q - 1)$-arc or a $(q - 2)$-arc. Let \\
$\theta^\prime_1$ be the number of planes $\Pi^\prime_i$, with $\Pi^\prime_i \not = ll^\prime$, containing a tangent of $\Pi_1 \cap K$ through $N_1$, a tangent of $\Pi_2 \cap K$ through $N_1$, where $\Pi^\prime_i \cap K$ is a $(q - 1)$-arc, \\
$\theta^\prime_2$ be the number of planes $\Pi^\prime_i$, with $\Pi^\prime_i \not = ll^\prime$, containing a tangent of $\Pi_1 \cap K$ through $N_1$, a tangent of $\Pi_2 \cap K$ through $N_1$, where $\Pi^\prime_i \cap K$ is a $(q - 2)$-arc. \\
Let $C^\prime_i \cap l^\prime = \lbrace P^\prime, N^\prime_i \rbrace$, with $l^\prime \cap K = \lbrace P^\prime \rbrace$ and $C^\prime_i$ the $(q + 2)$-arc extending $\Pi^\prime_i \cap K, i = 1, 2, ... , q + 1$. Then $N^\prime_i \not = N_1, i = 1, 2, ... , q + 1$. We may assume that $N^\prime_1 = N^\prime_2$. Assume, by way of contradiction, that $N^\prime_1 = N^\prime_2 = N^\prime_i$, with $i \in \lbrace 3, 4, ... , q + 1 \rbrace$. Then $N^\prime_1$ is on at least $3(q - 3) + 1$ tangents of $K$. So 
\begin{equation}
\label{eq41}
k \leq q^2 + q + 2 - 3q + 8 = q^2 -2q + 10.
\end{equation}
Hence
\begin{equation}
\label{eq42}
q^2 + (1 - \sqrt5)q + 5 < q^2 - 2q +10,
\end{equation}
that is,
\begin{equation}
\label{eq43}
(3 - \sqrt5)q < 5,
\end{equation}
clearly a contradiction. Hence we may assume that $N^\prime_1 = N^\prime_2, N^\prime_3 = N^\prime_4, N^\prime_1 \not = N^\prime_3, N^\prime_i \not \in \lbrace N^\prime_1, N^\prime_3 \rbrace$ for $i = 5, 6, ... , q + 1$. At least $\theta^\prime_1 - 4$ of the arcs $\Pi^\prime_5 \cap K, \Pi^\prime_6 \cap K, ... , \Pi^\prime_{q + 1} \cap K$ are $(q - 1)$-arcs, say $\Pi^\prime_5 \cap K, \Pi^\prime_6 \cap K, ... , \Pi^\prime_{\theta^\prime_1} \cap K$ are $(q - 1)$-arcs. The number of tangents of $\Pi^\prime_i \cap K$ containing $N^\prime_j$ , with $j \in \lbrace 1, 3 \rbrace$, is either 1 or 3, with $i = 5, 6, ... , \theta^\prime_1$; if $ N^\prime_j$ is contained in one tangent of $\Pi^\prime_i \cap K$, then $N^\prime_u$ is contained in 3 tangents of $\Pi^\prime_i \cap K$, with $\lbrace j, u \rbrace = \lbrace 1, 3 \rbrace$ and $i \in \lbrace 5, 6, ... , \theta^\prime_1 \rbrace$. So we may assume that at least $(\theta^\prime_1 - 4)/2$ of the $(q - 1)$-arcs $\Pi^\prime_i \cap K, i = 5, 6, ... , \theta^\prime_1$, have 3 tangents containing $N^\prime_1$. Counting the tangents of $K$ through $N^\prime_1$, we obtain at least
\begin{equation}
\label{eq44}
1 + (\theta^\prime_1 - 4) + (\theta^\prime_2 - 2) + 2(q - 3)
\end{equation}
tangents. As $\theta^\prime_1 + \theta^\prime_2 \geq q - 6$, this number of tangents is at least $1 + q - 6 - 6 + 2q - 6 = 3q - 17$. Hence 
\begin{equation}
\label{eq45}
k \leq q^2 + q + 2 - 3q + 17 = q^2 - 2q + 19.
\end{equation}
So
\begin{equation}
\label{eq46}
q^2 + (1 - \sqrt5)q + 5 < q^2 - 2q + 19,
\end{equation}
or
\begin{equation}
\label{eq47}
(3 - \sqrt5)q < 14,
\end{equation}
a contradiction for $q > 16$.

If at least one of the arcs $\Pi^\prime_1 \cap K, \Pi^\prime_2 \cap K$ is a $m$-arc with $m > q - 2$, then 
(\ref{eq44}) becomes
\begin{equation}
\label{eq48}
1 + (\theta^\prime_1 - 4) + (\theta^\prime_2 -1) + (q - 3) + (q - 2),
\end{equation}
which is at least $3q - 15$. Hence $k \leq q^2 - 2q + 17$. For $q = 16$ this gives $k \leq 241$. But for $q = 16$ the inequality $k > q^2 + (1 - \sqrt5)q + 5$ yields $k \ge 242$, a contradiction.

Finally we assume that $\Pi^\prime_1 \cap K$ and $\Pi^\prime_2 \cap K$ are $(q - 2)$-arcs. Then at least $\theta^\prime_1 - 2$ of the arcs $\Pi^\prime_i \cap K$, with $i = 5, 6, ... , q + 1$, are $(q - 1)$-arcs, say $\Pi^\prime_5 \cap K, \Pi^\prime_6 \cap K, ... , \Pi^\prime_{\theta^\prime_1 + 2} \cap K$. So at least $(\theta^\prime_1 - 2)/2$ of the $(q - 1)$-arcs $\Pi^\prime_i \cap K$, with $i = 5, 6, ... , \theta^\prime_1 + 2$, have 3 tangents containing either $N^\prime_1$ or $N^\prime_3$. First, assume that this is the case for $N^\prime_3$. If at least one of the arcs $\Pi^\prime_3 \cap K, \Pi^\prime_4 \cap K$ is a $m$-arc with $m > q - 2$, then the number of tangents of $K$ containing $N^\prime_3$ is at least
\begin{equation}
\label{eq49}
1 + (\theta^\prime_1 - 2) + (\theta^\prime_2 - 1) + (q - 3) + (q - 2),
\end{equation}
which is at least $3q - 13$. Hence $k \leq q^2 -2q + 15$, and so $q^2 + (1 - \sqrt5)q + 5 < q^2 - 2q + 15$, that is, $(3 - \sqrt5)q < 10$, a contradiction. Hence the arcs $\Pi^\prime_3 \cap K$ and $\Pi^\prime_4 \cap K$ are $(q - 2)$-arcs. Then the number of tangents of $K$ containing $N^\prime_3$ is at least
\begin{equation}
\label{eq50}
1 + \theta^\prime_1 + (\theta^\prime_2 - 2) + 2(q - 3),
\end{equation}
which is at least $3q - 13$. This yields again a contradiction. Consequently at least $(\theta^\prime_1 - 2)/2$ of the $(q - 1)$-arcs $\Pi^\prime_i \cap K$, with $i = 5, 6, ... , \theta^\prime_1 + 2$, have 3 tangents containing $N^\prime_1$. But then in (\ref{eq44}) $\theta^\prime_1 - 4$ may be replaced by $\theta^\prime_1 - 2$, yielding at least $3q - 15$ tangents of $K$ containing $N^\prime_1$. Hence $k \leq q^2 - 2q + 17$, which is a final contradiction.

We conclude that each tangent $l$ of $K$ is contained in at least one plane intersecting $K$ in at most four points.

(\rm C) {\bf Assume, by way of contradiction, that there is a tangent $l$ of $K$ which is contained in at least two planes $\Pi_1, \Pi_2$ intersecting $K$ in a $m$-arc, with $m \leq 4$}. \\

Assume that $l \cap K = \lbrace P \rbrace$ and that $\Pi_1 \cup \Pi_2$ contains $2q + \delta$ tangents of $K$ through $P$. We have $-5 \leq \delta \leq 1$.

(\rm C.1) {\bf Here we will show that $2q + \delta$ is the total number of tangents of $K$ containing $P$; as a corollary it will follow that $k \in \lbrace q^2 - q + 1, q^2 - q + 2, q^2 - q +3 \rbrace$}. Assume, by way of contradiction, that there is a tangent $l^\prime$ of $K$ containing $P$ with $l^\prime  \not \subset \Pi_1 \cup \Pi_2$. If $\vert ll^\prime \cap K \vert \leq 4$, then the number of tangents of $K$ containing $P$ is at least $2q + \delta + q - 3 = 3q + \delta - 3 \geq 3q - 8$, so $k \leq q^2 + q + 2 - 3q + 8 = q^2 - 2q + 10$. Hence
\begin{equation}
\label{eq51}
q^2 + (1 - \sqrt5)q + 5 < q^2 - 2q + 10,
\end{equation}
or $(3 - \sqrt5)q < 5$, a contradiction. Now we consider all planes containing $l^\prime$. By (B) at least one of these planes intersects $K$ in a $m$-arc, with $m \leq 4$. If at least two planes containing $l^\prime$ intersect $K$ in at most 4 points, then $P$ is contained in at least $2q - 5 + 2(q - 5) + 1 = 4q - 14$ tangents of $K$. Hence $k \leq q^2 + q + 2 - 4q + 14 = q^2 - 3q + 16$, so
\begin{equation}
\label{eq52}
q^2 + (1 - \sqrt5)q + 5 < q^2 - 3q + 16,
\end{equation}
that is, $(4 - \sqrt5)q < 11$, clearly a contradiction. Consequently exactly one plane $\Pi^\prime$ containing $l^\prime$ intersects $K$ in at most 4 points. Now we count the points of $K$ in the planes containing $l^\prime$. Let\\
$\theta_1$ be the number of planes , distinct from $ll^\prime$ and $\Pi^\prime$, containing $l^\prime$, containing a tangent of $K$ in $\Pi_1$ and containing a tangent of $K$ in $\Pi_2$,\\
$\theta_2$ be the number of planes containing $l^\prime$, distinct from $\Pi^\prime$, containing a tangent of $K$ in $\Pi_1$ and containing no tangent of $K$ in $\Pi_2$,\\
$\theta_3$ be the number of planes containing $l^\prime$, distinct from $\Pi^\prime$, containing a tangent of $K$ in $\Pi_2$ and containing no tangent of $K$ in $\Pi_1$,\\
$\theta_4$ be the number of planes , distinct from $\Pi^\prime$, containing $l^\prime$ and containing no tangent of $K$ in $\Pi_1$ or $\Pi_2$.\\
Then
\begin{equation}
\label{eq53}
k \leq 1 + (q - 1) + \theta_1(q - 2) + \theta_2(q - 1) + \theta_3(q - 1) + \theta_4q + 3,
\end{equation}
with
\begin{equation}
\label{eq54}
\theta_1 + \theta_2 + \theta_3 + \theta_4 = q - 1 \rm{and}\  \theta_2 + \theta_3 + 2\theta_4 \leq 6.
\end{equation}
So
\begin{equation}
\label{eq55}
k \leq q + 3 + (q - 1 - \theta_2 - \theta_3 - \theta_4)(q - 2) + \theta_2(q - 1) + \theta_3(q - 1) + \theta_4q,
\end{equation}
that is,
\begin{equation}
\label{eq56}
k \leq q^2 - 2q + 5 + (\theta_2 + \theta_3 + 2\theta_4),
\end{equation}
hence
\begin{equation}
\label{eq57}
k \leq q^2 - 2q + 5 + 6 = q^2 - 2q + 11.
\end{equation}
Consequently
\begin{equation}
\label{eq58}
q^2 + (1 - \sqrt5)q + 5 < q^2 - 2q + 11,
\end{equation}
or $(3 - \sqrt5)q < 6$, a contradiction.

It follows that $2q + \delta$ is the total number of tangents of $K$ containing $P$ and so $k = q^2 + q + 2 - 2q- \delta = q^2 - q + 2 - \delta$. As $k \leq q^2 - q + 3$ by Theorem 2.1, we have $-1 \leq \delta \leq 1$.

(\rm C.2) {\bf Some further properties of $K$}.
Let $l^{\prime\prime}$ be any tangent of $K$ not containing $P$ and let $K \cap l^{\prime\prime} = \lbrace P^\prime \rbrace$. By (B) $l^{\prime\prime}$ is contained in a plane $\Pi^{\prime\prime}$ with $\vert \Pi^{\prime\prime} \cap K \vert \leq 4$. There is a tangent $n$ of $K$ at $P^\prime$ not contained in $\Pi^{\prime\prime}$. The tangent $n$ is contained in a plane $\rho$ with $\vert \rho \cap K \vert \leq 4$. Let $2q + \delta^\prime$ be the number of tangents of $K$ at $P^\prime$ in $\rho \cup \Pi^{\prime\prime}$. Then $\delta^\prime \leq \delta$ and if $\rho \cap \Pi^{\prime\prime}$ is a tangent, then by the foregoing section we have $\delta^\prime = \delta$. Assume, by way of contradiction, that $\rho \cap \Pi^{\prime\prime}$ is not a tangent of $K$ and that $\delta^\prime < \delta$. Then there is a tangent $n^\prime$ of $K$ at $P^\prime$ not contained in $\rho \cup \Pi^{\prime\prime}$. The tangent $n^\prime$ is contained in a plane $\rho^\prime$ with $\vert \rho^\prime \cap K \vert \leq 4$. If $\rho \cap \rho^\prime$ is a tangent of $K$, then the $2q + \delta$ tangents of $K$ at $P^\prime$ are contained in $\rho \cup \rho^\prime$, a contradiction. So $\rho \cap \rho^\prime$ is not a tangent; similarly $\rho^\prime \cap \Pi^{\prime \prime}$ is not a tangent. Hence the number of tangents of $K$ at $P^\prime$ is at least $3(q - 2)$, so $2q + \delta \geq 3q - 6$, hence $\delta \geq q - 6$, a contradiction. We conclude that $\delta^\prime = \delta$ and that all tangents of $K$ at $P^\prime$ are contained in $\rho \cup \Pi^{\prime\prime}$.

Hence, given any point $Q \in K$ there are two planes $\alpha_1$ and $\alpha_2$ containing all tangents of $K$ at $Q$; also $\vert \alpha_1 \cap K \vert \leq 4$ and $\vert \alpha_2 \cap K \vert \leq 4$. These planes are uniquely defined by $Q$, and so is $\alpha_1 \cap \alpha_2$. By Section (A) the line $\alpha_1 \cap \alpha_2$ is a tangent of $K$ at $Q$. Let  $\widetilde \Pi$ be any plane containing $Q$, with $\widetilde \Pi \not \in \lbrace \alpha_1, \alpha_2 \rbrace$. Then $\widetilde \Pi \cap K$ contains at most two tangents at $Q$, so $\vert \widetilde \Pi \cap K \vert \geq q$. It follows that $K$ contains no $(q - 2)$-arcs and no $(q - 1)$-arcs.

Notice that $\vert \alpha_1 \cap K \vert + \vert \alpha_2 \cap K \vert + \delta = 3$ and remind that $-1 \leq \delta \leq1$.

Let $\widetilde \Pi$ be a plane containing $Q$, with $\widetilde \Pi \not \in \lbrace \alpha_1, \alpha_2 \rbrace$. The arc $\widetilde \Pi \cap K$contains always at least one tangent of $K$ at $Q$, except when $\delta = -1, k = q^2 - q + 3, \vert \alpha_1 \cap K \vert = \vert \alpha_2 \cap K \vert = 2$. So if $k \in \lbrace q^2 - q + 1, q^2 - q + 2 \rbrace$ and if $k = q^2 - q + 3$ with $\vert \alpha_1 \cap K \vert = \vert \alpha_2 \cap K \vert + 2 = 3$ or $\vert \alpha_2 \cap K \vert = \vert \alpha_1 \cap K \vert + 2 = 3$, then $\widetilde \Pi \cap K$ is not a $(q + 2)$-arc. If $\vert \alpha_1 \cap K \vert = \vert \alpha_2 \cap K \vert = 2, k = q^2 - q + 3$, then there is excactly one plane $\widetilde \Pi$ containing $Q$ for which $\widetilde \Pi \cap K$ is a $(q + 2)$-arc.\\

(\rm C.3)  \underline {$k = q^2 - q + 1$}\\
Then $\delta = 1$ and $\vert \Pi_1 \cap K \vert = \vert \Pi_2 \cap K \vert = 1$. Let $U_1, U_2 \in K$, with $U_1 \not = U_2$, and let $\xi_1, \xi_2$ be the planes containing $U_1$ intersecting $K$ in at most 4 points. If $U_2 \in \xi_1 \cup \xi_2$, then $\delta \leq 0$, a contradiction. Hence $U_2 \not \in \xi_1 \cup \xi_2$. Consequently any plane containing the line $U_1U_2$ has more than 4 points in common with $K$.

Now we count the points of $K$ in planes containing the line $U_1U_2$. Let $\theta_1$ be the number of planes containing $U_1U_2$ intersecting $K$ in a $q$-arc, and let $\theta_2$ be the number of planes containing $U_1U_2$ intersecting $K$ in a $(q + 1)$-arc. Then
\begin{equation}
\label{eq59}
\theta_1(q - 2) + \theta_2(q - 1) + 2 = q^2 - q + 1, {\rm with} \  \theta_1 + \theta_2 = q + 1.
\end{equation}
So
\begin{equation}
\label{eq60}
\theta_1(q - 2) + (q + 1 - \theta_1)(q - 1) + 2 = q^2 - q + 1,
\end{equation}
that is $\theta_1 = q$ and $\theta_2 = 1$.

Now we count the number of $(q + 1)$-arcs on $K$, and obtain
\begin{equation}
\label{eq61}
\frac{(q^2 - q + 1)(q^2 - q)}{(q + 1)q}.
\end{equation}
So $q + 1 \vert (q^2 - q + 1)(q - 1)$, so $q + 1 \vert 6$, a contradiction.

(\rm C.4) \underline {$k = q^2 - q + 2$}\\
Then $\delta = 0$ and $\lbrace \vert \Pi_1 \cap K \vert, \vert \Pi_2 \cap K \vert \rbrace = \lbrace 1, 2 \rbrace$. Let $Q$ be any point of $K$ and let $l_Q$ be the tangent of $K$ which is the intersection of the two planes $\alpha_1$ and $\alpha_2$ containing the $2q$ tangents of $K$ at $Q$. Let $(\alpha_1 \cup \alpha_2) \cap K = \lbrace Q, Q^\prime \rbrace$. Starting with $Q^\prime$ and $l_{Q^\prime}$, we find the same pair $\lbrace Q^\prime, Q \rbrace$. It follows that K is partitioned into pairs of type $\lbrace Q, Q^\prime \rbrace$. Let $\mathcal{L}$ be the set of these $(q^2 - q + 2)/2$ pairs.

Let $\lbrace Q, Q^\prime \rbrace \in \mathcal{L}$, let $\alpha_1$ and $\alpha_2$ be the planes containing the $2q$ tangents of $K$ at $Q$, and assume that $Q^\prime \in \alpha_1$. Then $\alpha_1 = l_Ql_{Q^\prime}$. Let $\Pi$ be a plane containing $QQ^\prime$, distinct from $\alpha_1$. As $\Pi$ contains a tangent of $K$ at $Q$, we have $\vert \Pi \cap K \vert \leq q + 1$. Counting the points of $K$ in the planes containing $QQ^\prime$, we obtain $\vert \Pi \cap K \vert = q + 1$. By an easy counting we see that the planes containing $l_Q$, but distinct from $\alpha_1$ and $\alpha_2$, intersect $K$ in $(q + 1)$-arcs. This way there arise $q - 1$ $(q + 1)$-arcs $K_1, K_2, ... , K_{q - 1}$, having kernels $N_1, N_2, ... , N_{q - 1}$ on $l_Q \setminus \lbrace Q \rbrace$. Assume, by way of contradiction, that $N_i = N_j, i \not = j$ and $i, j \in \lbrace 1, 2, ... , q - 1 \rbrace$. Then $N_i$ is on at least $2q + 1$ tangents of $K$, hence $k \leq q^2 - q + 1$, a contradiction. Let $l_Q \setminus \lbrace N_1, N_2, ... , N_{q - 1} \rbrace = N_Q$.

Assume, by way of contradiction, that $l_Q \cap l_{Q^\prime} \not = N_Q$. Let $l_Q \cap l_{Q^\prime} = N_i, i \in \lbrace 1, 2, ... , q - 1 \rbrace$, and let $R \in K_i \setminus \lbrace Q \rbrace$. Then $l_{Q^\prime} R \cap K$ is a $(q + 1)$-arc with kernel $N_i$. Hence $N_i$ is on at least $q^2 + 2$ tangents, a contradiction. Consequently $l_Q \cap l_{Q^\prime} = N_Q$; similarly, $l_Q \cap l_{Q^\prime} = N_{Q^\prime}$.

Assume, by way of contradiction, that $l_Q \cap l_S \not = \emptyset$, with $Q \not = S$ and $\lbrace Q, S \rbrace \not \in \mathcal{L}$. Let $\lbrace Q, Q^\prime \rbrace$ and $\lbrace S, S^\prime \rbrace$ be elements of $\mathcal{L}$. Now we count the number of tangents of $K$ containing $l_Q \cap l_S = M$. The arc $l_Ql_S \cap K$ is a $(q + 1)$-arc with kernel $M$, so $l_Ql_S$ contains $q + 1$ tangents of $K$ through $M$; the arc $l_QS^\prime \cap K$ is a $(q + 1)$-arc, and as the line $MS^\prime$ of the plane $l_SS^\prime$ is a tangent of $K$, the point $M$ is the kernel of $l_QS^\prime \cap K$, so $l_QS^\prime$ contains $q + 1$ tangents of $K$ through $M$; similarly the plane $l_SQ^\prime$ contains $q + 1$ tangents of $K$ through $M$. Hence $M$ is contained in more than $2q$ tangents of $K$, clearly a contradiction. It follows that if $l_Q \cap l_S \not = \emptyset$, with $Q \not = S$, then $\lbrace Q, S \rbrace \in \mathcal{L}$.

Let $\lbrace Q, S \rbrace \not \in \mathcal{L}$, with $Q$ and $S$ distinct points of $K$. Then $l_Q \cap l_S = \emptyset$. Now we count the points of $K$ in the planes containing the line $QS$. Let $\theta_1$ be the number of planes which contain $QS$ and intersect $K$ in a $q$-arc, and let $\theta_2$ be the number of planes which contain $QS$ and intersect $K$ in a $(q + 1)$-arc. Hence
\begin{equation}
\label{eq62}
\theta_1(q - 2) + \theta_2(q - 1) + 2 = q^2 - q + 2, {\rm with} \  \theta_1 + \theta_2 = q + 1.
\end{equation}
So $\theta_1 = q - 1$ and $\theta_2 = 2$. The 2 planes containing $QS$ and intersecting $K$ in a $(q + 1)$-arc are the planes $l_QS$ and $l_SQ$.

Let $\lbrace Q, S \rbrace \in \mathcal{L}$ and let $l_Q \cap l_S = N$. Then $N$ is kernel of no one of the $q - 1$ $(q + 1)$-arcs defined by planes containing the tangent $l_Q$ and of no one of the $q - 1$ $(q + 1)$-arcs defined by planes containing the tangent $l_S$. So for any line $n \not \in \lbrace l_Q, l_S \rbrace$ containing $N$ we have $\vert n \cap K \vert \in \lbrace 0, 2 \rbrace$. Let $n \cap K = \lbrace U, U^\prime \rbrace$.

First, assume that $\lbrace U, U^\prime \rbrace \not \in \mathcal{L}$. Then $\vert l_UU^\prime \cap K \vert = \vert l_{U^\prime} U \cap K \vert = q + 1$. As $\vert l_QU \cap K \vert = \vert l_SU \cap K \vert = q + 1$, the planes $l_QU$ and $l_SU$ are the two planes containing $UU^\prime$ and intersecting $K$ in a $(q + 1)$-arc. Hence $\lbrace l_QU, l_SU \rbrace = \lbrace l_UU^\prime, l_{U^\prime}U \rbrace$. So we may assume that $l_QU = l_UU^\prime$ and $l_SU = l_{U^\prime} U$. Consequently $l_Q \cap l_U \not = \emptyset$ and $l_S \cap l_{U^\prime} \not = \emptyset$, that is, $\lbrace Q, U \rbrace \in \mathcal{L}$ and $\lbrace S, U^\prime \rbrace \in \mathcal{L}$. Hence $\vert l_Ql_U \cap K \vert = \vert l_Sl_{U^\prime} \cap K \vert = 2$, clearly a contradiction as $Q, U, U^\prime \in l_Ql_U$.

It follows that $\lbrace U, U^\prime \rbrace \in \mathcal{L}$. So for any pair $\lbrace T, T^\prime \rbrace \in \mathcal{L}$, with $\lbrace T, T^\prime \rbrace \not = \lbrace Q, S \rbrace$, we have $N \in TT^\prime$. Let $n^\prime, n^{\prime \prime}$ be distinct lines containing $N$ with $n^\prime \not = n \not = n^{\prime \prime}$ and $n^\prime, n^{\prime \prime} \not \in \lbrace l_Q, l_S \rbrace$. Assume also that $n^\prime \cap K = \lbrace V, V^\prime \rbrace$ and $n^{\prime \prime} \cap K = \lbrace W, W^\prime \rbrace$. Then $\lbrace V, V^\prime \rbrace \in \mathcal{L}$ and $\lbrace W, W^\prime \rbrace \in \mathcal{L}$. By the foregoing the lines $VV^\prime, WW^\prime, QS$ contain $N$, clearly a contradiction.

(\rm C.5) \underline {$k = q^2 - q + 3$} \\
Let $P$ be any point of $K$ and let $l_P$ be the tangent of $K$ which is the intersection of the two planes $\Pi_1, \Pi_2$ containing the $2q - 1$ tangents of $K$ at $P$. Two cases are considered.

(\rm C.5.1) \underline {$\Pi_1 \cap K = \lbrace P, P^\prime, P^{\prime \prime} \rbrace, \Pi_2 \cap K = \lbrace P \rbrace$} \\
Then $K$ contains no plane $(q + 2)$-arcs containing $P$. Let $l$ be a tangent of $K$ at $P$, with $l$ in $\Pi_1$ and $l \not = l_P$. We count the points of $K$ in planes containing $l$. Let $\theta_1$ be the number of planes containing $l$ and intersecting $K$ in a $(q + 1)$-arc, and let $\theta_2$ be the number of planes containing $l$ and intersecting $K$ in a $q$-arc. Then 
\begin{equation}
\label{eq63}
\theta_1q + \theta_2(q - 1) + 3 = q^2 - q + 3, {\rm with} \ \theta_1 + \theta_2 = q.
\end{equation}
Hence $\theta_1 = 0$ and $\theta_2 = q$. Let $\widetilde{\Pi}_1, \widetilde{\Pi}_2, ... , \widetilde {\Pi}_q$ be the planes containing $l$ and intersecting $K$ in a $q$-arc, let $\widetilde{\Pi}_i \cap K = K_i$, let $C_i$ be the $(q + 2)$-arc extending $K_i$ and let $C_i \cap l = \lbrace P, N_i \rbrace$, with $i = 1, 2, ... , q$. Assume that for some $i \in \lbrace 1, 2, ... , q \rbrace$ we have $N_i \notin P^\prime P^{\prime \prime}$. The number of tangents of $K$ containing $N_i$ is at least 
\begin{equation}
\label{eq64}
q + (q - 1) + 2 = 2q + 1,
\end{equation}
a contradiction. Hence $N_1 = N_2 = \cdots = N_q = l \cap P^\prime P^{\prime \prime}$. Then the number of tangents of $K$ containing $N_1$ is at least
\begin{equation}
\label{eq65}
q(q - 1) + 1 = q^2 - q + 1,
\end{equation}
again a contradiction.

(\rm C.5.2) \underline {$\Pi_1 \cap K = \lbrace P, P^\prime \rbrace, \Pi_2 \cap K = \lbrace P, P^{\prime \prime} \rbrace$} \\
By (C.5.1), for each point $Q \in K$ the two planes $\alpha_1, \alpha_2$ through $Q$ intersecting $K$ in at most four points, intersect $K$ in exactly two points. If $\alpha_1 \cap K = \lbrace Q, Q^\prime \rbrace$ and $\alpha_2 \cap K = \lbrace Q, Q^{\prime \prime} \rbrace$, then the plane $QQ^\prime Q^{\prime \prime}$ is the only plane on $Q$ intersecting $K$ in a $(q + 2)$-arc. Hence the $(q + 2)$-arcs on $K$ partition $K$. So 
\begin{equation}
\label{eq66}
q + 2 \vert q^2 - q + 3, {\rm so} \ q + 2 \vert q - 7, {\rm so} \ q + 2 \vert 9,
\end{equation}
a contradiction.

Now the theorem is proved. \eop \\

\section{Corollaries} 

We are grateful to T. Sz\H{o}nyi for bringing reference \cite {SS: 93} to our attention which, in combination with Theorem 1.6, gives the following considerable improvement of the bound in Theorem 1.6; see also Remark 4.4.

\begin{theorem}
\label{thm 4.1} 
\begin{equation}
\label{eq67}
m^\prime_2(3, q) < q^2 - 2q +3\sqrt{q} + 2, q\ \text{even}, q \ge 2048.
\end{equation}
\end{theorem}
{\em Proof}. \quad
In \cite{SS: 93} it is proved that there does not exist a complete $k$-cap in $\PG(3, q)$, $q$ even, $q \ge 64$, for which
\begin{equation}
\label{eq68}
k \in [q^2 - (a - 1)q + a\sqrt{q} + 2 - a + \frac{a(a - 1)}{2}, q^2 - (a - 2)q - a^2\sqrt{q} ] 
\end{equation}
where $a$ is an integer which satisfies
\begin{equation}
\label{eq69}
2 \le a \le \frac{-2\sqrt{q} + 3 + \sqrt{16q\sqrt{q} + 12q - 44\sqrt{q} - 7}}{4\sqrt{q} + 2}.
\end{equation}
Putting $a = 3$, the desired result immediately follows from Theorem 1.6. \eop\\

\begin{theorem}
\label{thm 4.2}
\begin{itemize}
\item [{\rm(i)}] $m_2(4, 4) = 41$,
\item [{\rm(ii)}] $m_2(4, 8) \leq 479$,
\item [{\rm(iii)}] $m_2(4, q) < q^3 - q^2 + 2\sqrt{5}q - 8$, $q$ even, $q > 8$.
\end{itemize}
\end{theorem}
{\em Proof}. \quad
For $q = 4$, see \cite {YE: 99}. Assume, by way of contradiction, that $K$ is a $k$-cap of $\PG(4, 8)$ with $k > 479$, or a $k$-cap of $\PG(4, q)$, $q$ even and $q > 8$, with $k > q^3 - q^2 + 2\sqrt{5}q - 8$. At each of its points the cap $K$ has $t = q^3 + q^2 + q + 2 - k$ tangents. Hence we assume that $t < 107$ for $q = 8$ and $t < 2q^2 + (1 - 2\sqrt5)q + 10$ for $q > 8$. We obtain a contradiction in several stages.

{\rm I} {\bf $K$ contains no plane $q$-arc} \\
Similar to the reasoning in Section I in the proof of Theorem 6.27 in \cite {JWPH: 16}.

{\rm II} {\bf There exists no solid $\delta$ such that $q^2 + 1 > \vert \delta \cap K \vert > q^2 + (1 - \sqrt5)q + 5$} \\
Suppose $\delta$ exists. Let $\delta \cap K = K'$. Then $K'$ can be completed to an ovoid $O$ of $\delta$, by Theorem 3.1. Let $N \in O \setminus K'$ and let $N' \in K'$. Consider the $q + 1$ planes of $\delta$ through $NN'$. Since each of these planes meets $O$ in a $(q + 1)$-arc, each plane meets $K'$ in at most a $q$-arc. By I, there is no $q$-arc on $K$; so each plane meets $K'$ in at most a $(q -1)$-arc.

Assume, by way of contradiction, that none of these intersections is a $(q - 1)$-arc. Therefore a count of the points on $K'$ gives
\begin{equation}
\label{eq70}
\vert K' \vert \leq (q + 1)(q - 3) +1,
\end{equation}
whence
\begin{equation}
\label{eq71}
q^2 + (1 - \sqrt5)q + 5 < q^2 - 2q - 2,
\end{equation}
so
\begin{equation}
\label{eq72}
(3 - \sqrt5)q + 7 < 0,
\end{equation}
a contradiction.

So we may assume that for one of the planes $\delta$ through $NN'$, say $\Pi$, we have $\vert \Pi \cap K' \vert = q - 1$. Now we consider all solids of $\PG(4, q)$ containing the plane $\Pi$. Let $\theta$ be the number of solids $\Pi'$ for which $\vert \Pi' \cap K \vert > q^2 + (1 - \sqrt5)q + 5$, so $q + 1 - \theta$ is the number of solids $\Pi''$ for which $\vert \Pi'' \cap K \vert < q^2 + (1 - \sqrt5)q + 5$. We have $\theta \geq 1$.

First, assume $\theta \geq 2$. So there are at least two solids $\Pi'_1, \Pi'_2$ containing $\Pi$ such that $\vert \Pi'_i \cap K \vert > q^2 + (1 - \sqrt5)q + 5$, with $i = 1, 2$. By Theorem 3.1 $\Pi'_i \cap K$ can be completed to an ovoid $O_i$ of $\Pi^\prime_i, i = 1, 2$. So $O_i \cap \Pi$ is a $(q + 1)$-arc $(\Pi \cap K') \cup \lbrace N'_i, N''_i \rbrace, i = 1, 2.$ Since $\Pi \cap K'$ can be contained in no more than three $(q + 1)$-arcs, contained in a common $(q + 2)$-arc, we have $\vert \lbrace N'_1, N''_1 \rbrace \cap \lbrace N'_2, N''_2 \rbrace \vert \geq 1$. Assume $N'_1 = N'_2$. So the number of tangents of $K$ containing $N'_1$ is at least 
\begin{equation}
\label{eq73}
2(q^2 + (1 - \sqrt5)q + 5 - q + 1) + q - 1,
\end{equation}
so
\begin{equation}
\label{eq74}
2q^2 + (1 - 2\sqrt5)q + 11,
\end{equation}
a contradiction.

Finally, assume that $\theta = 1$. Counting the points of $K$ in the $q + 1$ solids, we obtain
\begin{equation}
\label{eq75}
k < q(q^2 + (1 - \sqrt5)q + 5 - q + 1) + (q^2 - 1),
\end{equation}
that is,
\begin{equation}
\label{eq76}
k < q^3 + (1 - \sqrt5)q^2 + 6q - 1.
\end{equation}
Hence, for $q > 8$,
\begin{equation}
\label{eq77}
q^3 - q^2 + 2\sqrt{5}q - 8 < q^3 + (1 - \sqrt5)q^2 + 6q - 1,
\end{equation}
so
\begin{equation}
\label{eq78}
0 < (2 - \sqrt5)q^2 + (6 - 2\sqrt5)q + 7,
\end{equation}
a contradiction. For $q = 8$, there arises $479 < 479$, a contradiction.

{\rm III} {\bf For a point $N$ not in $K$, there do not exist planes $\Pi_1$ and $\Pi_2$ such that $\Pi_1 \cap \Pi_2 = N$ and such that $\Pi_i \cap K$ is a $(q + 1)$-arc with nucleus $N$ for $i = 1, 2$} \\
Similar to the reasoning in Section III in the proof of Theorem 6.27 in \cite {JWPH: 16}.

{\rm IV} {\bf The tangents through any point $Q$ off $K$ lie in a solid} \\
Similar to the reasoning in Section IV in the proof of Theorem 6.27 in \cite {JWPH: 16}.

{\rm V} {\bf The final contradiction is obtained by counting the tangents of $K$} \\
Similar to the reasoning in Section V in the proof of Theorem 6.27 in \cite {JWPH: 16}. \eop \\

\begin{theorem}
\label{thm 4.3}
For $q$ even, $q > 2$, $n \geq 5$,
\begin{itemize}
\item [{\rm(i)}] $m_2(n, 4) \leq \frac{118}{3}.4^{n - 4} + \frac{5}{3}$
\item [{\rm(ii)}] $m_2(n, 8) \leq 478.8^{n - 4} - 2(8^{n - 5} +\cdots+ 8 + 1) + 1$,
\item [{\rm(iii)}] $m_2(n, q) < q^{n - 1} - q^{n - 2} + 2\sqrt{5}q^{n - 3} - 9q^{n - 4} - 2(q^{n - 5} +\cdots+ q + 1) + 1$, for $q > 8$.
\end{itemize}
\end{theorem}
{\em Proof} \quad
This follows directly from Theorem 1.1, Theorem 4.2 and Theorem 6.14(ii) in \cite {JWPH: 16}. \eop \\

\begin{remark}
\label{rem 4.4}
{\rm The bound in Theorem 4.1 leads to considerable improvements of Theorem 4.2 and Theorem 4.3. We just   mention these bounds, but the proofs are the theme of a subsequent paper}. \\

{\rm For $q$ even, $q \ge 2048$},

\begin{equation}
\label{eq79}
m_2(4, q) < q^3 - 2q^2 + 3q \sqrt{q} + 8q - 9\sqrt{q} - 6.
\end{equation}
{\rm For $q$ even, $q \ge 2048$, $n \ge 5$},

\begin{equation}
\label{eq80}
m_2(n, q) < q^{n - 1} - 2q^{n - 2} + 3q^{n -3}\sqrt{q} + 8q^{n - 3} - 9q^{n - 4}\sqrt{q} - 7q^{n -4} - 2(q^{n - 5} + \cdots + q + 1) + 1. 
\end{equation}

\end{remark}

\section{Remark} 
The bound in the MAIN THEOREM is better than the bound of Chao, see \cite {JMC: 99}. In 2014 Cao and Ou, see \cite {JMC: 14}, published the bound $k < q^2 - 2q + 8$ ($q$ even and $q \geq 128$), which is better than ours. I did not follow some reasoning in their proof, so I sent two mails to one of the authors explaining why I think Section 1.3 of the proof is not correct. Unfortunately I never received an answer. 




\end{document}